\documentclass[11pt]{article}

\usepackage{amssymb}
\usepackage{amsfonts}

\title{Sprague-Grundy theory in bounded arithmetic\\(Preliminary Draft)}
\author{Satoru Kuroda\\
Gunma Prefectural Women's University}
\date{}

\newtheorem{definition}{Definition}
\newtheorem{theorem}{Theorem}
\newtheorem{lemma}{Lemma}
\newtheorem{proposition}{Proposition}
\newtheorem{corollary}{Corollary}

\newcommand{\ltwo}{\mathcal{L}_A^2}
\newcommand{\lpg}{\mathcal{L}_{NK}}

\newcommand{\tpg}{V_{NK}}

\newcommand{\sbi}[1]{\Sigma^B_{#1}}

\newcommand{\all}{\forall}
\newcommand{\ext}{\exists}
\newcommand{\pairs}[1]{\langle #1\rangle}

\begin{document}

\maketitle

\begin{abstract}
In this paper, we formalize Sprague-Grundy theory for combinatorial games 
in bounded arithmetic. We show that in the presence of Sprague-Grundy numbers, 
a fairly weak axioms capture PSPACE. 
\end{abstract}

\section{Introduction}

Since the seminal paper by Bouton \cite{bouton}, combinatorial games have been paid  
much attention in various branches of mathematics. The observation in \cite{bouton} is 
later generalized by Grundy \cite{grundy} and Sprague \cite{sprague} to form 
a powerful tool for finding winning strategies which is called Grundy number or 
Sprague-Grundy number. 

Deciding the complexity of perfect information games is also a major problem 
in computational complexity theory. Many combinatorial games are related to 
space complexity such as PSPACE. For instance, Schaefer \cite{schaefer} proved that 
the game Node Kayles played on undirected graphs is complete for PSPACE. 
while some games have much weaker complexity such as P or LOGSPACE. 

In this paper we show that with the aid of Sprague-Grundy number, a fairly weak theory 
of two-sort bounded arithmetic can capture PSPACE. More precisely, we introduce 
a function computing Sprague-Grundy number for Node Kayles together with strategy functions 
for both players using Sprague-Grundy number to the system $V^0$ 
and show that any alternating polynomial time machine can be simulated by a game of Node Kayles. 

Specifically, for an alternating Turing machine $M$ and an input $X$, we construct in $V^0$ 
an undirected graph $G(M,X)$ such that Alice has an winning strategy if and only if 
$M$ accepts $X$. Since the strategy functions are polynomial time computable in Sprague-Grundy 
function, this result suggests that Sprague-Grundy number has such a strong 
computational power that manages search through polynomial space. 

There are a number of literature concerning bounded arithmetic for PSPACE.  
Buss \cite{buss} in his seminal paper defined a second order theory $U^1_2$ 
whose provably total functions coincide with PSPACE. Later, Skelley \cite{skelley} 
defined a three sort system $W^1_1$ for PSPACE. While these theories require higher order 
objects compared to theories for classes inside the polynomial hierarchy, Eguchi \cite{eguchi} 
defined a PSPACE theory $\sbi{0}$-ID by extending the two sort language by predicates 
which represent inductive definition for $\sbi{0}$ definable relations. 
Our theory $\tpg$ presented in this paper is considered as a minimal theory for PSPACE 
as it is contained in any of the above theory. We also remark that an application of bounded 
arithmetic to combinatorial game theory is also given by Soltys and Wilson \cite{sw} 
who showed that strategy stealing argument can be formalized in $W^1_1$ and in turn proved that 
the game Chomp is in PSPACE. 

We can alternatively formalize our theory with a stronger base theory such as $PV$ while 
introducing Sprague-Grundy function only. However we do not follow such an approach since 
formalizing in weak theory such as $V^0$ enables us to construct theories for combinatorial 
games having weaker computational power. Among such games we are particularly interested in 
the game NIM whose computational complexity is around LOGSPACE but no completeness result 
is known so far. We remark that this choice of base theory forces us to give a slightly 
more complicated construction of the graph $G(M,X)$.

There is a rich theory of combinatorial games with a number of games and so we hope that 
our result gives a neat framework for logical analysis of combinatorial games. 

This paper is organized as follows: in section 2 we define our theory $\tpg$ by extending $V^0$ 
by functions computing winning strategies. In section 3, we show that $\tpg$ actually computes 
winning strategies for Node Kayles. Section 4 is devoted to the proof of our main theorem. 
In particular, we construct a graph so that players winning strategies witness accepting or 
rejecting computations. 

\section{Formalizing combinatorial games}

We will formalize the argument for combinatorial games in the language of 
two-sort bounded arithmetic.

We will assume familiarity with basic notions and properties of two-sort 
bounded arithmetic. For a detail, readers should consult with textbooks such as 
\cite{cn}. 

Let $\ltwo$ be the two sort language of Cook-Nguyen \cite{cn}.  
Basically, upper case letters denote binary strings and lower case letters denote 
natural numbers. We also adopt an unusual notation that vector presentation of 
lower case letters such as $\bar{z}$ also denote strings. 
For a language $L$ we denote the $\sbi{0}$ formulas in $L$ by $\sbi{0}(L)$. 

The theory $V^0$ has defining axioms for symbols in $\ltwo$ together with the 
bit-comprehension axiom for $\sbi{0}$ formulas. We use many properties of $V^0$ 
in this paper whose details can be found in \cite{cn}. 

For a string $X$ and a number $i<|X|$, $X(i)$ denotes both the predicate that the $i$th 
bit of $i$ is $1$ and the $i$th bit of $X$ itself. The sequence of numbers are coded by 
a string and we define the $i$th entry of a sequence $X$ by $X[k]$ and the length of $X$ by 
$Len(X)$. For two sequences $X$ and $Y$, we denote the concatenation by $X*Y$. 
Strings are sometimes identified with a binary sequence as 
$P=\pairs{p_0,\ldots,p_n}$. Coding such sequences and proving basic properties of 
sequences can be done in $V^0$.

The game we consider is known as Node Kayles which is played over undirected graphs. 
We code graphs by a two-dimensional array where we assume that any node has an edge to 
itself. Two-dimensional arrays represent directed graphs in general and undirected graphs 
are given as a symmetric relation which is coded by symmetric matrices. So we define
$$
\begin{array}{l}
DGraph(G)\Leftrightarrow\\
\all x\in G\ext u,v<|G|(x=\pairs{u,v})\land\\
(\all u<|G|(\ext v<|G|(\pairs{u,v}\in G\lor\pairs{v,u}\in G))\rightarrow\pairs{u,u}\in G).\\
UGraph(G)\Leftrightarrow DGraph(G)\land\all u,v<|G|(\pairs{u,v}\in G\rightarrow\pairs{v,u}\in G).\\
Node(G)=V_G=\{u<|G| \:\ \pairs{u,u}\in G\}.
\end{array}
$$

We define the game Node-Kayles over undirected graphs to be 
an impartial game played by two players Alice and Bob (Alice always moves first) 
starting from a graph $G$ and in the move with the option $G'$ which is a subgraph of $G$, 
the player chooses a node $x\in Node(G)$ and returns the subgraph $G_x$ which is 
defined by 
$$
Node(G'_x)=\{y\in Node(G') : \pairs{x,y}\not\in E_{G'}\}
$$
and
$$
\pairs{y,z}\in G'_x\Leftrightarrow y,z\in Node(G'_x)\land\pairs{y,z}\in G'.
$$
For a sequence $\bar{w}=\pairs{w_1,\ldots,w_l}$ we define $G_{\bar{w}}$ inductively as 
$$
G_\emptyset=G,\ 
G_{\bar{w}*v}=\left\{
\begin{array}{ll}
(G_{\bar{w}})_v&\mbox{ if }v\in Node(G_{\bar{w}})\\
G_{\bar{w}}&\mbox{ otherwise.}
\end{array}
\right.
$$

The first player unable to move loses. So a game over $G$ is coded by a sequence 
$\bar{w}=\pairs{w_1,\ldots,w_l}$ such that $w_1\in G$, $w_{i+1}\in G_{\pairs{w_1,\ldots,w_i}}$ 
for any $i<l$, $G_{w_{l-1}}\neq\emptyset$ and $G_{w_l}=\emptyset$. 
Alice wins in the game $W$ over $G$ if $Len(w)\bmod 2= 1$ and otherwise Bob wins. 

\begin{proposition}
The function computing $G_{\bar{w}}$ from $G$ and $\bar{w}$ is $\sbi{0}$-definable in $V^0$. 
\end{proposition}

\noindent
(Proof). It is easy to see that $G_{\bar{w}}$ is definable by the formula
$$
\begin{array}{l}
\varphi(G,G',\bar{w})\Leftrightarrow\\
\all u,v<|V_G|(\pairs{u,v}\in G'\leftrightarrow
(\pairs{u,v}\in G\land\neg\ext w_i\in\bar{w}(w=u)).
\end{array}
$$
so that $\all G,\bar{w}\ext! G'\varphi(G,G',\bar{w})$ is provable in $V^0$. \hfill$\Box$.

Now we will define our base theory for combinatorial games. 
First we introduce functions $sg(G)$, $\tau(G)$, $\tau_A(\pairs{b_0,\ldots,b_l},G)$ and 
$\tau_B(\pairs{a_0,\ldots,a_l},G)$ with the following defining axioms:
$$
\begin{array}{l}
G=\emptyset\rightarrow sg(G)=0,\\
\neg UGraph(G)\rightarrow sg(G)=\max\{x\in V_G\}+1,\\
UGraph(G)\land G\neq\emptyset\rightarrow
sg(G)=\min\{k<|V_G|:\all x\in V_Gk\neq sg(G_x)\}.
\end{array}
$$
$$
\tau(G)=
\left\{
\begin{array}{ll}
\min\{v\in V_G\ :\ sg(G_v)=0\}&\mbox{ if such }v\mbox{ exists.}\\
\max\{v\in V_G\}+1&\mbox{ otherwise.}
\end{array}
\right.
$$
$$
\begin{array}{l}
\tau_A(\emptyset,G)=\tau(G),\\
\tau_A(\pairs{b_0,\ldots,b_{l+1}},G)=
\tau_A(\pairs{b_0,\ldots,b_l},G)*\pairs{b_{l+1},\tau(G_{\tau_A(\pairs{b_0,\ldots,b_l},G)*b_{l+1}})}
\end{array}
$$
$$
\begin{array}{l}
\tau_B(\emptyset,G)=\emptyset,\\
\tau_B(\pairs{a_0,\ldots,a_{l+1}},G)=
\tau_B(\pairs{a_0,\ldots,a_l},G)*\pairs{a_{l+1},\tau(G_{\tau_B(\pairs{a_0,\ldots,a_l},G)*a_{l+1}})}.
\end{array}
$$

\begin{definition}
Let $\lpg$ be the language $L_A^2$ extended by function symbols $sg(G)$ 
$\tau(G)$, $\tau_A(\pairs{b_0,\ldots,b_l},G)$ and $\tau_B(\pairs{a_0,\ldots,a_l},G)$ . 
The $\lpg$ theory $\tpg$ comprises the following axioms:
\begin{itemize}
\item defining axioms for symbols in $\lpg$ 
\item $\sbi{0}(\lpg)$-COMP: $\ext X<a\all y<a(X(a)\leftrightarrow\varphi(a))$, \\
where $\varphi(a)\in\sbi{0}(\lpg)$ which does not contain free occurrences of $X$. 
\end{itemize}
\end{definition}

Thus $\tpg$ is $V^0$ in the extended language $\lpg$. 

\noindent
{\bf Remark.} We need only functions $sg$ and $\tau_A$ in order to axiomatize 
the theory $\tpg$ since other two functions are definable from these functions. 
For instance, $\tau_G$ can be defined from $sg$ and $\tau_B$ can be defined by $\tau_A$. 
However we add these two functions to the language to make argument simple. 

The following fact is well-known.
\begin{proposition}
$\tpg$ proves $\sbi{0}(\lpg)$-IND:
$$
\varphi(0)\land
\all x(\varphi(x)\rightarrow\varphi(x+1))\rightarrow\all x\varphi(x).
$$
\end{proposition}






\section{Winning strategies in Sprague-Grundy system}

We show that strategy functions $\tau_A$ and $\tau_B$ actually computes winning game instances 
for Alice and Bob respectively. 

\begin{definition}
Define formulas $AWS_{\tau_A}(G,l)$ and $BWS_{\tau_A}(G,l)$ as follows:
$$
\begin{array}{ll}
AWS_{\tau_A}(G)\Leftrightarrow&
\all l\all\pairs{b_0,\ldots,b_l}[(l=\lfloor|V_G|/2\rfloor\land\all i<l b_i\leq|V_G|+1)\\
&\phantom{\all l\all\pairs{b_0,\ldots,b_l}[}
\rightarrow\ext l_0\leq l(\all i<l_0(b_i\in Node(G_{\tau_A(\pairs{b_0,\ldots,b_{i-1}},G)})\land\\
&\phantom{\all l\all\pairs{b_0,\ldots,b_l}[\rightarrow}
\tau(G_{\tau_A(\pairs{b_0,\ldots,b_{i-1}},G)})\in Node(G_{\tau_A(\pairs{b_0,\ldots,b_{i-1}},G)}))\\
&\phantom{\all l\all\pairs{b_0,\ldots,b_l}[\rightarrow}
\land b_{l_0}\not\in Node(G_{\tau_A(\pairs{b_0,\ldots,b_{l_0-1}},G)}))]
\end{array}
$$
$$
\begin{array}{ll}BWS_{\tau_B}(G)\Leftrightarrow&
\all l\all\pairs{a_0,\ldots,a_l}[(l=\lceil|V_G|/2\rceil\land\all i<l a_i\leq|V_G|+1)\\
&\phantom{\all l\all\pairs{a_0,\ldots,a_l}[}
\rightarrow\ext l_0\leq l(\all i<l_0(a_i\in Node(G_{\tau_A(\pairs{a_0,\ldots,a_{i-1}},G)})\land\\
&\phantom{\all l\all\pairs{b_0,\ldots,b_l}[\rightarrow}
\tau(G_{\tau_A(\pairs{a_0,\ldots,a_{i-1}},G)})\in Node(G_{\tau_A(\pairs{a_0,\ldots,a_{i-1}},G)}))\\
&\phantom{\all l\all\pairs{a_0,\ldots,a_l}[\rightarrow}
\land a_{l_0}\not\in Node(G_{\tau_A(\pairs{a_0,\ldots,a_{l_0-1}},G)}))]
\end{array}
$$
\end{definition}

\begin{theorem}
$\tpg$ proves that 
$$
\all G\left\{
Ugraph(G)\rightarrow((sg(G)\neq 0\rightarrow AWS_{\tau_A}(G)\land 
(sg(G)=0\rightarrow BWS_{\tau_B}(G))))\right\}.
$$
\end{theorem}

\noindent
(Proof). 
We argue inside $\tpg$. 

Suppose that $sg(G)\neq 0$ and let $\pairs{b_0,\ldots,b_l}$ be a list of nodes 
in $G$ where $l=\lfloor|V_G|/2\rfloor$. We show that 
$$
\all i<l(\all j\leq i b_j\in Node(G_{\tau_A(\pairs{b_0,\ldots,b_{j-1}},G)})
\rightarrow sg(G_{\tau_A(\pairs{b_0,\ldots,b_i},G)}=0)\mbox{ (*)} 
$$
The proof proceeds by induction on $i$. 

If $i=0$ then (*) trivially follows by the assumption. 

Suppose by the inductive hypothesis that (*) holds for $i\geq 0$ and assume that 
$$
b_{i+1}\in Node(G_{\tau_A(\pairs{b_0,\ldots,b_i},G)}).
$$
Since $sg(G_{\tau_A(\pairs{b_0,\ldots,b_i},G)})=0$, it must be that 
$$
sg(G_{\tau_A(\pairs{b_0,\ldots,b_i},G)}*b_{i+1})\neq 0
$$
and by the definition of $\tau$, we have
$$
sg(G_{\tau_A(\pairs{b_0,\ldots,b_{i+1}},G)})=0.
$$
So we have (*) for $i+1$. 

We argue similarly for the case of $sg(G)=0$ and by noting that (*) is a $\sbi{0}$ 
formula,  the claim is obtained by $\sbi{0}$-IND in $\tpg$. \hfill$\Box$

\section{Sprague-Grundy system captures PSPACE}
\label{sec:sprag-grundy-syst}

Now we are ready to show our main result; the theory $\tpg$ captures $PSPACE$. 

\begin{theorem}
A function is $\sbi{1}$ definable in $\tpg$ if and only if it is in PSPACE.
\end{theorem}

\noindent
(Proof). It is easy to show that functions $sg$, $\tau$, $\tau_A$ and $\tau_B$ can be 
computed in PSPACE. So the only if part can be proved using the standard witnessing 
argument. Actually the provably total functions of the universal conservative extension 
of $\tpg$ is the $AC^0$ closure of functions $sg$ and $\tau_A$. So Herbrand theorem 
implies the witnessing. 
Thus the proof of if part is given is the rest of this section. \hfill$\Box$

We will show that any polynomial time alternating Turing machine can be simulated 
by a game in $\tpg$. First recall that PSPACE is equal to APTIME 
(cf. Papadimitriou \cite{papa}).  So we actually show 
that any polynomial-time alternating Turing machine can be simulated by a game of 
Node Kayles. 

We assume some harmless simplifications on alternating Turing machines. 
Let $M$ be an alternating Turing machine with time bound $p(|X|)$ on input $X$. 
where we assume that $p(n)$ is even for all $n$. We assume that 
all computation of $M$ on input $X$ terminates exactly at time $p(|X|)$. 
We also assume that the space bound of $M$ is $p(|X|)$. 
Furthermore, we assume that $M$ is binary branching. So we formalize 
the transition function as 
$$
\delta_M(k,q,a)=\pairs{q_k,a_k,m_k}
$$
where $k=0,1$, $q$ and $q_k$ are states of $Q$ and $a,a_k\leq 2$, $m_k\in\{-1,0,1\}$. 
We abuse the notation and write 
$$
\delta_M(k,C,C')\Leftrightarrow
C'\mbox{ is the next configuration of }C\mbox{ along the path }k.
$$
The final assumption is that $M$ computes in normal form in the sense that 
it first guesses the path $P=\pairs{p_1,\ldots,pi_{p(n)}}$ in the computation 
tree and then start computing using $P$. 

We show that polynomial time bounded alternating Turing machines can be simulated 
by Node Kayles provably in $\tpg$. 

Let $C_{INIT}(M,X)$ denote the initial configuration of $M$ on input $X$. 
For a binary string $P$, we denote by $C(P,M,X)$ the configuration of $M$ reachable from 
$C_{INIT}(M,X)$ along the path $P$. The predicate $Accept(C,M)$ denotes that 
$C$ is an accepting configuration of $M$. Note that all these functions 
and predicates are definable in $V^0$. We also define
$$
\begin{array}{l}
Comp(\pairs{C_0,\ldots,C_{p(|X|)}},P,M,X)\Leftrightarrow\\
\phantom{Comp}C_0=C_{INIT}(M,X)\land\all i<p(|X|)\delta_M((P)_i,C_i,C_{i+1}), \\
Acomp(\pairs{C_0,\ldots,C_{p(|X|)}},P,M,X)\Leftrightarrow\\
\phantom{Comp}Comp(\pairs{C_0,\ldots,C_{p(|X|)}},P,M,X)\land Accept(C_{p(|X|)},M,X),\\
Rcomp(\pairs{C_0,\ldots,C_{p(|X|)}},P,M,X)\Leftrightarrow\\
\phantom{Comp}Comp(\pairs{C_0,\ldots,C_{p(|X|)}},P,M,X)\land\neg Accept(C_{p(|X|)},M,X),\\
\end{array}
$$

\begin{theorem}\label{theorem:main}
There exist functions $G(M,X)$, $Comp_A(M,X,P)$, $Comp_R(M,X,P)$, $Path_A(M,X,P)$ and 
$Path_R(M,X,P)$ which are $\sbi{1}$ definable in $\tpg$ such that the following formulas are 
provable in $\tpg$. 
$$
\begin{array}{l}
\mbox{(1).}\ \all M,X UGraph(G(M,X)),\\
\mbox{(2).}\ \all M,X,P(|P|=p(|X|)/2\rightarrow(Len(Path_A(M,X,P))=2Len(P)\land\\
\phantom{\mbox{(2).}\ }\all k<Len(Path_A(M,X,P))(Path_A(M,X,P)[2k+1]=P[k])),\\
\mbox{(3).}\ \all M,X,P(|P|=p(|X|)/2\rightarrow(Len(Path_R(M,X,P))=2Len(P)\land\\
\phantom{\mbox{(2).}\ }\all k<Len(Path_R(M,X,P))(Path_R(M,X,P)[2k]=P[k])),\\
\mbox{(4).}\ 
\all M,X\\
\phantom{\mbox{(2).}\ }\{(sg(G(M,X))\neq 0\rightarrow\all P(|P|=p(|X|)/2\rightarrow \\
\phantom{\mbox{(2).}\ }AComp(Comp_A(M,X,Path_A(M,X,P)),Path_A(G(M,X),P),M,X)))\\
\phantom{\mbox{(2).}\ }\land(sg(G(M,X))=0\rightarrow\all P(|P|=p(|X|)/2\rightarrow \\
\phantom{\mbox{(2).}\ }RComp(Comp_R(M,X,Path_R(M,X,P)),Path_R(G(M,X),P),M,X)))\}\\
\end{array}
$$
\end{theorem}

First we sketch the outline of the proof. 

Let $M$ be an alternating Turing machine and $X$ be an input. We construct two graphs 
$G_A(M,X)$ and $G_B(M,X)$ so that each legitimate game instance of either games corresponds to 
a computation of $M$ on input $X$. 
Specifically, the first $p(|X|)$ moves of the game constitute a path $P$ with $|P|=p(|X|)$ 
followed by a list of moves which establishes a computation of $M$ along the path $P$, 
if players move correctly. We require that $G_A(M,X)$ and $G_R(M,X)$ satisfy that 
a game instance $I$ is A-winning if and only if $I$ corresponds to an accepting and rejecting 
computation of $M$ on $X$ along $P$ respectively. 

Once the graph is constructed, we can extract functions $Comp_A(G,P)$, $Comp_B(G,P)$, 
$Path_A(G,P)$ and $Path_B(G,P)$ using strategy functions $\tau_A$ and $\tau_B$. 

Now we present details of the proof. 

The construction of $G_A(M,X)$ and $G_R(M,X)$ is similar to that for the graph 
simulating QBF games in \cite{schaefer}. 
Let $M=(Q,\Sigma,\delta,q_0,q_A)$ be an alternating Turing machine with 
$Q=\{q_0,\ldots,q_m\}$, $\Sigma=\{0,1,2\}$ where $2$ denotes the blank symbol and $q_A=q_1$. 
The transition function is given as $\delta(p,q,a)=\pairs{q_p,a_p,m_p}$ where $p\in\{0,1\}$, 
$q,q_p\in Q$ and $m_p\in\{-1,0,1\}$ whose intended meaning is that if the current state is $q$, 
the head reads the symbol $a$ and the path $p$ is chosen then the state changes to $q_P$, 
the tape content of the current head position is overwritten by $a_P$ and the head moves 
by $m_P$.  

Let $s=p(|X|)$ be the number of alternations of $M$ on $X$, $l_0=p(|X|)+2$ be the length 
of the sequence coding configurations and $n_0=2(s+1)l_0$. It turns out that $s+n_0$ is 
equal to the number of total moves in A-winning legitimate game instances.  
We construct the graph  of $G_A(M,X)$ and $G_R(M,X)$ with layers $P_i,A_{i,j},B_{i,j},$ 
of legitimate nodes, $Y_i$ of illegitimate nodes and $C_A,C_B$ of constraints nodes 
so that in the $i$th round, the player must choose her or his move from $i$th legitimate 
layer. Nodes in each layers are given as follows:

\begin{itemize}
\item $P$-layers $P_i=\{p_{i,0},p_{i,1}\}$ for $0\leq i<s$ represent the choice of $i$th 
path in the computation. 
\item $A$-layers $A_{i,j}$ corresponds to computation by Alice after the path is decided by 
choices from $P_0,\ldots,P_{s-1}$ and consists of nodes as follows:
$$
\begin{array}{l}
A_{i,j}=\{a_{i,j,0}^T,a_{i,j,1}^T,a_{i,j,2}^T\},\ 0\leq j<s\\
\vspace{3pt}
A_{i,s}=\{a_{i,k}^H\ :\ 0\leq k<s\}\\
\vspace{3pt}
A_{i,s+1}=\{a_{i,r}^Q\ :\ 0\leq r<|Q|\}\\
\end{array}
$$
The intended meaning is that if Alice chooses nodes $a_{i,0,i_0}^T,\ldots,a_{i,s-1,i_{s-1}}^T$, 
$a_{i,k}^H$ and $a_{i,r}^Q$ then Alice's computation of the $i$th configuration is 
$C_i=\pairs{q_r,k,i_0,\dots,i_{s-1}}$. 
\item $B$-layers $B_{i,j}=\{b_{i,j}\}$ which are intended for Bob's moves for 
$0\leq i<s$ and $0\leq j\leq s+1$ or $i=s$ and $0\leq j\leq s$. Note that 
Bob's have no choice of moves for these rounds. Also note that the number of $B$-layers 
is one less than that of $A$-layers. 
\end{itemize}

We list these layers in the order that players choose their moves as 
$$
P_0,\ldots,P_{s-1},A_{0,0},B_{0,0},\ldots,A_{s,s+1},B_{s,s}.
$$
So we sometimes denote layers by ignoring their types as 
$$
L_k=\left\{
\begin{array}{ll}
P_k&\mbox{ if }0\leq k<s,\\
A_{i,j}&\mbox{ if }k=s+2(i\cdot l_0+j),\ 0\leq i\leq s,\ 0\leq j\leq s+1,\\
B_{i,j}&\mbox{ if }k=s+2(i\cdot l_0+j)+1,\ 0\leq i\leq s,\ 0\leq j\leq s.\\
\end{array}
\right.
$$

We define constraint layers $C_A$ and $C_R$ for $G_A(M,X)$ and $G_R(M,X)$ respectively 
which expresses constraints for the computation of $M$. 
Nodes of these layers are labelled by propositional formulas and we identify nodes with 
their labels. The layer $C_A$ and $C_R$ contain the following nodes: 
\begin{description}
\item[(A)] Nodes of the first sort are called initial nodes and express 
the initial configuration of $M$ on $X$ which consists of $\rightarrow a_{0,0}^Q$, 
$\rightarrow a_{0,0}^H$, for $j<|X|$, $\rightarrow a_{0,j,k}^T$ where $k=X(j)$ and for 
$|X|\leq j<s$, $\rightarrow a_{0,j,2}^T$. 
\item[(B)] The second sort are called transition nodes of $M$ which consists of rules expressing 
the transition function of $M$. 
Specifically, let $c\in\{0,1\}$, $0\leq j\leq m$, $z\in\{0,1,2\}$ and 
$\delta(c,q_j,z)=\pairs{q_{j'},z',d}$ for some $0\leq j\leq|Q|$, 
$z'\in\{0,1,2\}$ and $d\in\{-1,0,1\}$. 
Then for $0\leq i<s$ $0\leq j\leq|Q|$ and $0\leq k<s$ , we introduce the following rules:
$$
\begin{array}{l}
p_{i,c}\land a_{i,j}^Q\land a_{i,k}^H\land a_{i,k,z}^T\rightarrow a_{i+1,k,z'}^T\\
p_{i,c}\land a_{i,j}^Q\land a_{i,k}^H\land a_{i,k',a}^T\rightarrow a_{i+1,k',a}^T,\ k'\neq k\\
p_{i,c}\land z_{i,j}^Q\land z_{i,k}^H\land z_{i,k,a}^T\rightarrow a_{i+1,k+d}^H\\
p_{i,c}\land a_{i,j}^Q\land a_{i,k}^H\land a_{i,k,a}^T\rightarrow a_{i+1,j'}^Q\\
\end{array}
$$
Note that these rules compute the $i+1$st configuration from the $i$th configuration 
which is specified by choosing the path $c$. 
We call a rule containing $p_{i,c}$ for $c=0,1$ as $i$-rule. 
\end{description}
Moreover, $C_A$ contains a single accepting node denoted by $Acc$ while $C_R$ contains 
a single rejecting node denoted by $Rej$. 

Finally, the non-legitimate nodes are defined as 
$$
Y_{n_0-k}=\{y_{n_0-k,n_0-k+j}\ :\ 0\leq j<k+1\}.
$$
for $1\leq k<n_0$. 

Next we define edges among the nodes.  In the following, let $C$ denote either $C_A$ or $C_R$. 
\begin{enumerate}
\item For $0\leq i<s$ and $c\in\{0,1\}$, $p_{i,c}\in P_i$ is connected to all nodes 
in $C$ which contains $p_{i,1-c}$. 
\item For $0\leq i\leq s$ and $0\leq j\leq s+1$, $a\in A_{i,j}$ is connected to all nodes in 
$C$ which either contain $a$ in the succedent or $b\in A_{i,j}$ with $b\neq a$ in the antecedent. 
\item The node $a_{s,1}^Q$ in $G_A(M,X)$ is connected to the node $Acc$. 
\item The node $a_{s,j}^Q$ for $j\neq 1$ in $G_A(M,X)$ is connected to the node $Rej$. 
\item all nodes in $C$ are mutually connected. 
\item All nodes in $L_k\cup Y_k$ for $1\leq k\leq t_0$ are mutually connected. 
\item The node $y_{t_0-k,t_0-k+j}\in Y_{t_0,k}$ is connected to all nodes in 
$$
\bigcup\{L_i\cup Y_i\ :\ t_0-k<i\leq t_0+1,\ i\neq t_0-k+j\}.
$$
\end{enumerate}

\begin{proposition}
The function computing $G_A(M,X)$ and $G_R(M,X)$ from $M$ and $X$ is $\sbi{1}$ definable in $V^0$. 
\end{proposition}

\noindent
(Proof). We code $G(M,X)$ in such a way that indices of nodes represent their labels. 
For instance, the node $p_{i,c}$ in $P_i$ for $0\leq i<s$ and $c\in\{0,1\}$ is indexed 
by the tuple $\pairs{0,i,c}$ where the first entry $0$ represents that it belongs to a $P$-layer. 

Similarly, the node $a_{i,j}^Q$ in $A_{i,s+1}$ for $0\leq i\leq s$ and $0\leq j\leq |Q|$ 
is indexed by the tuple $\pairs{0,s+i\cdot n_0+1,j}$ and nodes in other $A$-layers and 
$B$-layers are indexed as well. 

The node $y_{n_0-k,n_0-k+j}$ in the layer $Y_{n_0-k}$ is indexed by the tuple 
$\pairs{1,n_0-k,j}$ for $0\leq j<k+1$. 

Finally nodes in $C_A\cup C_R$ are indexed by tuples of the form 
$\pairs{0,n_0,t}$ where $t$ is a tuple coding its label. For instance the node 
$$
p_{i,c}\land a_{i,j}^Q\land a_{i,k}^H\land a_{i,h,a}^T\rightarrow a_{i+1,j_{c,a}}^Q
$$
is denoted by the tuple $\pairs{0,i,c,j,k,a,0,j_{c,a}}$.

Then it is easy to see that the edge relation of $G(M,X)$ is definable by a $\sbi{0}$ 
formula so it is defined by $\sbi{0}$-COMP.\hfill$\Box$

\vspace{6pt}

We say that a subgraph $G'$ of $G=G_A(M,X)$ or $G_R(M,X)$ is $k$-legitimate for 
$0\leq k\leq s+n_0$ if 
$$
Leg(G',G,k)\Leftrightarrow
\all x\in V_G
\left(
\left(
x\in\bigcup_{k'<k}L_k\rightarrow x\not\in V_{G'}
\right)
\land
\left(
x\in\bigcup_{k<k'\leq s+n_0}L_k\rightarrow x\in V_{G'}
\right)
\right).
$$

In the following, we denote $G=G_A$ or $G_R$ if there is no fear of confusion.

The following lemma states that the graph $G(M,X)$ is constructed 
so that players are forced to choose their moves from legitimate nodes for otherwise they lead to 
an immediate loose. 
\begin{lemma}\label{lemma:legitimate1}
$V_{NK}$ proves that from any legitimate graph $G'$ of $G$, the first non-legitimate move leads to 
an immediate lose for either player:
$$
\all G'\all <n_0+s\all x((Leg(G',G,k)\land x\not\in L_{k+1})\rightarrow sg(G'_x)\neq 0). 
$$ 
\end{lemma}

\noindent
(Proof). We argue in $\tpg$ to show that if $G'$ is a $k$-legitimate subgraph of $G(M,X)$ 
and $v\not\in L_k$ then $sg(G'_v)\neq 0$. 

Let $v\not\in L_k$. Then either $v\in Y_j$ for $j\geq k$ or $v\in L_j$ for $j>k$. 
In the first case, we have $v=y_{j,l}$ for some $l$ and taking it from $G'$ removes 
all nodes except $L_l\cup Y_l$. Since $L_l\cup Y_l$ forms a complete subgraph, 
it must be that $sg(G_{j,l})\neq 0$. 

In the second case, $G_v$ consists of all nodes in $L_l\cup Y_l$ with $l\neq j$ and 
nodes in $L_{N_0+1}$ which are not connected to $v$. By the construction of $G(M,X)$, 
$y_{k,j}$ remains in $G'_v$ and is connected to all nodes in $G'_v$. So we have 
$G'_{\pairs{v,y_{k,j}}}=\emptyset$. This implies that $sg(G'_v)\neq 0$ as required. \hfill$\Box$


\vspace{6pt}

We say that a sequence $\bar{w}=\pairs{v_1,\cdots,v_m}$ of nodes in $G_A(M,X)$ or 
$G_R(M,X)$ is legitimate, denoted by $SLeg(w,G)$, if $v_i\in L_i$ for all $i\leq m$. 
Then the following is an immediate consequence of Lemma \ref{lemma:legitimate1}. 

\begin{corollary}\label{cor:legitimate2}
$\tpg$ proves that 
$$
\all k<t_0\all\pairs{v_1,\ldots,v_k}\mbox{ : legitimate }\all v_{k+1}
(v_{k+1}\not\in L_{k+1}\rightarrow sg(G_{v_1\cdots v_{k+1}})\neq 0).
$$
\end{corollary}

\noindent
(Proof). 
It remains to show that if $\pairs{v_1,\ldots,v_m}$ is legitimate then 
for any $k\leq m$, $G_{v_1\cdots v_k}$ is a $k$-legitimate subgraph of $G(M,X)$ 
which can be proved by $\sbi{0}$-IND on $k\leq m$. \hfill$\Box$


\vspace{6pt}

If both players move legitimately, The first $s$ moves will be $p_{0,c_0},\ldots,p_{s-1,c_{s-1}}$ 
which decides the path $P=\pairs{c_0,\ldots,c_{s-1}}$ in the computation tree of $M$ on $X$. 

We require that if $sg(G_A(M,X)_P)=0$ then Bob can win the game for $G_A(M,X)_P$ 
only if he moves consistently with the computation of $M$ on $X$ along the path $P$. 
Otherwise if $sg(G_A(M,X)_P)\neq 0$ then Alice can win the game for $G_R(M,X)_P$ only if 
she moves consistently with the computation along $P$. 

In order to prove the above property of $G(M,X)$ in $\tpg$, 
we next show that each list of legitimate moves forms a list of configurations. 

Note that we can divide A-layers and B-layers into consecutive lists $A_{i,0},\dots,A_{i,s+1}$ and 
$B_{i,0},\dots,B_{i,s+1}$. We call these two lists as the $i$-round. We assert 
that each set of legitimate move by both Alice and Bob for the $i$-round forms a configuration of $M$ 
on input $X$. Specifically, let Alice's moves for the $i$-round be given as 
$$
\bar{a}_i=a_{i,j}^Q, a_{i,k}^H, a_{i,0,a_0}^T,\ldots,a_{i,s-1,a_{s-1}}^T.
$$
Then we define $conf(\bar{a}_i)=\pairs{j,k,a_0,\ldots,a_{s-1}}$. 
Thus a legitimate sequence $\pairs{\bar{a}_0,\ldots,\bar{a}_s}$ of moves by Alice 
forms a sequence of configurations $\pairs{conf(\bar{a}_0),\ldots,conf(\bar{a}_s)}$. 

We define legitimate moves by Alice and Bob after $s+2$ rounds as 
$$
\begin{array}{l}
Leg(\pairs{v_1,\ldots,v_k},M,X)\Leftrightarrow\all j<k(v_{j+1}\in L_{s+1+i}),\\
A\mbox{-}Leg(\pairs{a_0,\ldots,a_k},M,X)\Leftrightarrow\all j<k(a_{j+1}\in L_{s+2j}),\\
B\mbox{-}Leg(\pairs{b_0,\ldots,b_k},M,X)\Leftrightarrow\all j<k(b_{j+1}\in L_{s+2j+1}).
\end{array}
$$
We omit parameters $M$ and $X$ if it is clear from the context. We also denote legitimate 
sequences of Alice and Bob as 
$\pairs{a_{0,0},a_{0,1},\ldots,a_{i,j}}\mbox{ and }\pairs{b_{0,0},b_{0,1},\ldots,b_{i,j}}$ 
respectively for $i\leq s$ and $j\leq s+1$. 

Finally we define predicates which states that a given legitimate move form 
a computation of $M$. 
$$
\begin{array}{l}
Comp(\pairs{a_{0,0},\ldots,a_{s,s+1}},M,X,P)
\Leftrightarrow\\
\phantom{Comp}Leg(\bar{a},M,X)\land conf(\bar{a}_0)=C_{INIT}(M,X)\land
\all i<s\delta_M(P(i),conf(\bar{a}_i),conf(\bar{a}_{i+1})),\\
AComp(\pairs{a_{0,0},\ldots,a_{s,s+1}},M,X,P)
\Leftrightarrow\\
\phantom{Comp}Comp(\pairs{a_{0,0},\ldots,a_{s,s+1}},M,X,P)\land Accept(\bar{a}_s,M,X),\\
RComp(\pairs{a_{0,0},\ldots,a_{s,s+1}},M,X,P)
\Leftrightarrow\\
\phantom{Comp}Comp(\pairs{a_{0,0},\ldots,a_{s,s+1}},M,X,P)\land\neg Accept(\bar{a}_s,M,X).
\end{array}
$$
Note that Bob's moves after $s$ rounds are unique if he moves legitimately. 
So we denote $\bar{b}=\pairs{b_{0,0},\ldots,b_{s,s}}$. 

In the followings, $M$ and $X$ always denote a code of an alternating TM and its 
input respectively and we refrain from stating it explicitly. 

For a sequence $X=\pairs{x_0,\ldots,x_l}$, We define the function 
$ASeq(X)=\{x_i:i\bmod 2=0\}$. Note that if $X$ codes a game instance then 
$ASeq(X)$ gives a list of Alice's moves. 

The next lemma states that the value of $sg(G_A(M,X)_P)$ for $|P|=s$ decides whether $M$ accepts $X$ 
along the path $P$. 
\begin{lemma}\label{lemma:comp2}
$\tpg$ proves that 
$$
\begin{array}{l}
\all M,X,P\biggl\{|P|=s\rightarrow\\
(sg(G_A(M,X)_P)\neq 0\rightarrow
AComp(ASeq(\tau_A(\pairs{b_{0,0},\ldots,b_{s,s}},G_A(M,X)_P),M,X,P)))\land\\
(sg(G_A(M,X)_P)=0\rightarrow
RComp(ASeq(\tau_A(\pairs{b_{0,0},\ldots,b_{s,s}},G_R(M,X)_P),M,X,P)))\biggl\}. 
\end{array}
$$
\end{lemma}

In order to prove Lemma \ref{lemma:comp2}, we first prepare some notations. 
As stated above, legitimate moves $\bar{a}_i$ by Alice in $a_i$ rounds is presented as
$$
\bar{a}_i=a_{i,0,k_0}^T,\ldots,a_{i,s-1,k_{s-1}}^T,a_{i,k}^H,a_{i,j}^Q
$$
where $0\leq j\leq m$, $0\leq k\leq s-1$ and $k_0\ldots,k_{s-1}\in\{0,1,2\}$. 
Likewise, Bob's moves for $a_i$ rounds is represented as
$\bar{b}_i=b_{i,1},b_{i,2},\ldots,b_{i,s+2}$ for $0\leq i<s$ and 
$\bar{b}_s=b_{s,1},b_{s,2},\ldots,b_{s,s+1}$. 

We denote the moves by Alice and Bob for $G(M,X)_P$ respectively as
$$
\bar{a}=\pairs{\bar{a}_0,\ldots, \bar{a}_s}\mbox{ and }
\bar{b}=\pairs{\bar{b}_0,\ldots, \bar{b}_s}
$$

We sometimes ignore the type of the nodes of Alice's move and denote by 
$a_{i,j}$ the $j$-th move of Alice in the $i$-round. Furthermore we define 
$$
\bar{a}^{\leq i,j}=\bar{a}_1\ldots,\bar{a}_{i-1},a_{i,1},\ldots,a_{i,j}\mbox{ and }
\bar{a}^{<i,j}=\bar{a}_1\ldots,\bar{a}_{i-1},a_{i,1},\ldots,a_{i,j-1}.
$$
$$
\bar{a}^{\leq i}=\bar{a}_1\ldots,\bar{a}_{i}\mbox{ and }
\bar{a}^{<i}=\bar{a}_1\ldots,\bar{a}_{i-1}.
$$
The sequences $\bar{b}^{\leq i,j}$, $\bar{b}^{<i,j}$, $\bar{b}^{\leq i}$ and 
$\bar{b}^{<i}$ are defined similarly. 

For sequences $\bar{a}=\pairs{a_0,\ldots,a_k}$ and $\bar{b}=\pairs{b_0,\ldots,b_k}$ 
or $\pairs{b_0,\ldots,b_{k-1}}$, we define the $V^0$-definable function  
$$
merge(\bar{a},\bar{b})=\pairs{a_0,b_0,\ldots,a_k,b_k}
\mbox{ or }\pairs{a_0,b_0,\ldots,a_k,b_k}.
$$
respectively.

The proof of Lemma \ref{lemma:comp2} is divided into a series of sublemmas.
Define $\sbi{0}$ formulas $Init(r,z,M,X)$ and $Next(r,z,p,C,M)$ so that 

$$
\begin{array}{l}
Init(r,z,M,X)\Leftrightarrow z\mbox{ is the }r\mbox{th element of }C_{INIT}(M,X),\\
Next(r,z,p,C,M)\Leftrightarrow z\mbox{ is the }r\mbox{th element of }C'\mbox{ with }\delta_M(p,C,C').\\
\end{array}
$$
A $A_{i,j}$-rule is a transition rule in $C_A$ whose succedent contains a node in $A_{i,j}$. 
We say that a legitimate subgraph $G'$ of $G(M,X)$ contains no $A_{i,j}$-rule 
if there is no node in $G'$ which belongs to $C_A$ and 
represents some $A_{i,j}$-rule. We also say that $G'$ contains no $A$-rules if for all $i\leq s$ 
and $j\leq s+1$, $G'$ contains no $A_{i,j}$-rules. 
Note that thess properties are formalized by a $\sbi{0}$ formula. 

Let $\bar{z}=\pairs{z_0,\ldots,z_k}$ be a list of legitimate moves by Alice or Bob for $k\leq (s+1)(s+2)$. 
We define that $\bar{z}$ is a partial computation as
$$
\begin{array}{l}
PComp(\bar{a},P,M,X)\Leftrightarrow\\
A\mbox{-}Leg(\bar{a},G_P)\land
\all k\leq Len(\bar{a})
\{(q_k=0\rightarrow Init(r_k,a_k,M,X))\land\\
(q_k>0\rightarrow Next(r_k,a_k,P(q_k-1),conf(\bar{a}_{q_k-1}),M)),
\end{array}
$$
where $q_k$ and $r_k$ are such that $k=q_k(s+1)+r_k$ and $0\leq r_k\leq s+1$. 

The next lemma states that  moves by Alice or Bob must be consistent with the computation of $M$ in order 
to obtain legitimate options. 
\begin{lemma}
Let $G$ be either $G_A(M,X)$ or $G_R(M,X)$. Then $\tpg$ proves that 
$$
\begin{array}{l}
\all M,X\all P\all l\leq (s+1)l_0
\all\bar{a}=\pairs{a_0,\ldots,a_l}\all\bar{b}=\pairs{b_0,\ldots,b_l}\\
\biggl\{(|P|=s\land A\mbox{-}Leg(\bar{a})\land B\mbox{-}Leg(\bar{b}))\land Len(\bar{a})=Len(\bar{b})+1
\rightarrow\\
(PComp(\bar{a},P,M,X)\leftrightarrow
\all k\leq l(G_{P*merge(\bar{a},\bar{b})}
\mbox{ contains no }A_{q_k,r_k}\mbox{-rule}))
\biggr\}. 
\end{array}
$$
\end{lemma}

\noindent
(Proof). We prove the claim of the lemma for $A_{i,j}$-rules by induction on $l$. 
If $l=0$ then we have to do nothing. So suppose that $l\geq 0$ and by the 
inductive hypothesis assume that the claim holds for $l$. 
Let us denote the lefthand side of the subformula inside the brace $\{\cdots\}$ of the claim 
by $(*)_l$. Assume that $(*)_{l+1}$ holds, 
that is 
$$
\all k\leq l+1((q_k=0\rightarrow Init(r_k,a_k,M,X))\land
(q_k>0\rightarrow\delta_M(p_{q_{k-1}},conf(\bar{a}_{q_{k-1}}),2r_k-1,a_k)).
$$
By the inductive hypothesis we already have 
$$
\all k\leq l((G_P)_{merge(\bar{a},\bar{b})}\mbox{ contains no }A_{q_k,r_k}\mbox{-rules}).
$$
So it suffice to show that 
$(G_P)_{merge(\bar{a},\bar{b})}$ contains no $A_{q_{l+1},r_{l+1}}$-rules 

If $q_{l+1}=0$ then we have $Init(r_{l+1},a_{l+1},M,X)$ and since $\rightarrow a_{l+1}$ 
is the only $L_{l+1}$-rule, we have the claim. Otherwise, we have
$$
Next(2r_{l+1}-1,a_{l+1},p_{q_{l+1}-1},conf(\bar{a}_{q_{l+1}-1}),M)
$$
so there must be a rule in $C$ of the form $A\rightarrow a_{l+1}$ where 
$A$ represents a conjunction which is consistent with $conf(\bar{a}_{q_{l+1}-1})$. 
Furthermore, it is the only $A_{q_{l+1},r_{l+1}}$-rule which is in 
$(G_P)_{merge(\bar{a}_{\leq l},\bar{b}_{\leq l})}$. 
Thus again we have the claim. 

Conversely, suppose that $(*)_{l+1}$ does not hold. If $(*)_l$ does not hold then 
we have the claim by the inductive hypothesis. So suppose that 
$$
\begin{array}{l}
(q_{l+1}=0\land\neg Init(r_{l+1},a_{l+1},M,X))\\
\lor(q_{l+1}>0\land\neg Next(2r_{l+1}-1,a_{l+1},p_{q_{l+1}-1},conf(\bar{a}_{q_{l+1}-1}),M)).
\end{array}
$$
If the first disjunct is true then there exists an initial rule $\rightarrow y_{l+1}$ 
where $y_{l+1}\in L_{l+1}$ and $y_{l+1}\neq a_{l+1}$ which is not eliminated by the move 
$a_{l+1}$ of Alice. 

Otherwise if the second conjunct is true then we may assume that $(G_P)_{merge(\bar{a},\bar{b})}$ 
does not contain any $L_k$-rule for $k\leq l$. Since 
$$
\neg Next(2r_{l+1},y_{l+1},p_{q_{l+1}-1},conf(\bar{a}_{q_{l+1}-1}),M)
$$
there must be a rule of the form $A\rightarrow y_{l+1}$ such that $A$ is consistent with 
$conf(\bar{a}_{q_{l+1}-1})$ and so it remains in $(G_P)_{merge(\bar{a},\bar{b})}$. 
Since $A\rightarrow y_{l+1}$ is not eliminated by $a_{l+1}$ we have the claim. 
\hfill$\Box$


\begin{corollary}\label{cor:accepting1}
Let $G$ be either $G_A(M,X)$ or $G_R(M,X)$. 
Then $\tpg$ proves that if Alice moves legitimately on $G_P$ then she removes all $A$-rules 
if and only if her moves are consistent with the computation of $M$ on $X$ along $P$:  
$$
\begin{array}{l}
\all M,X,P\all\bar{a}
\biggl\{(|P|=s\land Leg(\bar{a})\land Len(\bar{a})=n_0)\rightarrow\\
(Comp(\bar{a},P,M,X)
\leftrightarrow(G_{P*\pairs{e_0,e}*merge(\bar{a},\bar{b})}\mbox{ contains no $A$-rules of }M)
\biggr\}
\end{array}
$$
\end{corollary}

\noindent
(Proof). 
We argue inside $\tpg$. First we remark that 
\begin{itemize}
\item the move $a_{i,j}$ by Alice removes all nodes in $C$ which contain $a_{i,j}$ in the 
succedent or  $a'\in X_{i,j}$ with $a'\neq a_{i,j}$ in the antecedent and 
\item any move in $\bar{b}$ by Bob does not remove any node in $C$. 
\end{itemize}
We say that a node in $C$ is a $i$-round rule if it is a $A_{i,j}$-rule for some $0\leq j\leq s+1$. 
We will prove that 
$$
\begin{array}{l}
conf(\bar{a}_0)=C_{INIT}(M,X)\rightarrow
\all k\leq N\{\all i\leq k\delta_M(P(i),conf(\bar{a}_i),conf(\bar{a}_{i+1}))\\
\leftrightarrow
\all i\leq k(G_P)_{merge(\bar{a}^{\leq i},\bar{b}^{\leq i})}\mbox{ contains no $i$-round rules}\}. 
\end{array}
$$
The proof is by induction on $k$. For $k=0$ we show that 
$$
\begin{array}{l}
conf(\bar{a}_0)=C_{INIT}(M,X)\\
\leftrightarrow 
(G_P)_{merge(\bar{a}_0,\bar{b}_0)}\mbox{ contains no initial rule of }M.
\end{array}
$$
Suppose first that $conf(\bar{a}_0)=C_{INIT}(M,X)$. Then each move $a_{0,j}$ of Alice 
removes the initial rule $\rightarrow a_{0,j}$ in $L_{N_0}$. Such a rule exists since 
$conf(\bar{a}_0)=C_{INIT}(M,X)$. 

Conversely, suppose that $conf(\bar{a}_0)\neq C_{INIT}(M,X)$. Then for some choice 
$a_{0,j}$ of Alice, $L_{N_0}$ contains the initial rule $\rightarrow z'_{0,j}$ with 
$a_{0,j}\neq z'_{0,j}$. Since $\rightarrow z'_{0,j}$ cannot be removed by any other moves 
in $a_0$-rounds, $(G_P)_{merge(\bar{a}_0,\bar{b}_0)}$ must contain it. 

For induction step, suppose that for $k\leq s-1$
$$
\begin{array}{l}
(conf(\bar{a}_0)=C_{INIT}(M,X)\land
\all i<k\delta_M(P(i),conf(\bar{a}_i),conf(\bar{a}_{i+1}))\\
\leftrightarrow
\all i<k(G_P)_{merge(\bar{a}^{\leq k},\bar{b}^{\leq k})}\mbox{ contains no $i$-round rules}. 
\end{array}
$$
and we show that
$$
\delta_M(P(i),conf(\bar{a}_i),conf(\bar{a}_{i+1}))\leftrightarrow
(G_P)_{merge(\bar{a}^{\leq k+1},\bar{b}^{\leq k+1})}\mbox{ contains no $k$-round rules}. 
$$

Suppose that $\delta_M(P(i),conf(\bar{a}_i),conf(\bar{a}_{i+1}))$ holds. 
By the construction of $G(M,X)$, antecedents of $k+1$ rules of 
$(G_P)_{merge(\bar{a}^{\leq k},\bar{b}^{\leq k})}$ form $conf(\bar{a}_k)$. 

In $a_{k+1}$-rounds, Alice must choose nodes in order to remove all such nodes in $L_{N_0}$. 
Since each such node specifies a transition rule of $M$, we have the claim. 

Also the induction step is easily seen by the above remarks. Since the claim 
is $\sbi{0}$, it is proved by $\sbi{0}$-IND in $\tpg$ and the claim of the lemma 
easily immediately follows. \hfill$\Box$

\vspace{6pt}

Let $G$ be a graph and $z_0\ldots,z_k\in V_G$. We say that $\pairs{z_0,\ldots,z_k}$ is 
a winning sequent for $G$, denoted by $WSeq(\pairs{z_0,\ldots,z_k},G)$ if 
$$
G_{\pairs{z_0,\ldots,z_{k-1}}}\neq\emptyset\land G_{\pairs{z_0,\ldots,z_k}}=\emptyset.
$$

\begin{corollary}\label{cor:accepting2-1}
$\tpg$ proves that Alice's moves for $G_A(M,X)_P$ form an accepting computation 
if and only if Alice wins the game:
$$
\begin{array}{l}
\all M,X,P\all\bar{a}=\pairs{a_{0,0},\ldots a_{s,s+1}}
\biggl\{(|P|=s\land Leg(\bar{a})\land Len(\bar{a})=(s+1)(s+2))\rightarrow\\
(AComp(\bar{a},P,M,X)\leftrightarrow WSeq(merge(\bar{a},\bar{b}),G_A(M,X)_P)\biggr\}.
\end{array}
$$
\end{corollary}

\noindent
(Proof). First note that Bob cannot removes any nodes in $C_A$ unless he can move 
legitimately for a node in $C$. 
By Lemma \ref{cor:accepting1}, the only node in $C_A$ which may remain 
in $(G_P)_{merge(\bar{a},\bar{b})}$ is the acceptance node $Acc$. 
So we have 
$$
conf(\bar{a}_s)=C_{ACCEPT}(M,X)
\leftrightarrow Acc\mbox{ is removed in $a_s$-rounds}".
$$


\begin{corollary}\label{cor:accepting2-2}
$\tpg$ proves that Alice's moves for $G_R(M,X)_P$ form a rejecting computation 
if and only if Alice wins the game:
$$
\begin{array}{l}
\all M,X,P\all\bar{a}=\pairs{a_{0,0},\ldots a_{s,s+1}}
\biggl\{(|P|=s\land Leg(\bar{a})\land Len(\bar{a})=(s+1)(s+2))\rightarrow\\
(RComp(\bar{a},P,M,X)\leftrightarrow WSeq(merge(\bar{a},\bar{b}),G_R(M,X)_P)\biggr\}.
\end{array}
$$
\end{corollary}

\noindent
(Proof). The proof is almost identical to Corollary \ref{cor:accepting2-1}. 
The only difference is if Alice moves in accordance with the computation of $M$ 
on $X$ along $P$ then she must remove the rejecting node $Rej$ by the last move. 
\hfill$\Box$

\vspace{6pt}

In order to show that the strategy function yields computations of $M$, we need to 
relate Sprague-Grundy number of $G=G_A(M,X)$ or $G_R(M,X)$ and the computation of $M$. 
The next lemma asserts that Alice can always chooses options $G'$ of $G_A(M,X)_P$ 
so that $sg(G')=0$ if and only if Alice's moves form an accepting computation along $P$.
\begin{lemma}\label{lemma:accepting}
$\tpg$ proves that 
$$
\begin{array}{l}
\all M,X,P\all\bar{a}=\pairs{a_{0,0},\ldots a_{s,s+1}}
\biggl\{(|P|=s\land Leg(\bar{a}))\rightarrow\\
\all k<Len(\bar{a})(sg(G_{P*merge(\bar{a}^{\leq k},\bar{b}^{< k})})=0)\leftrightarrow
AComp(\bar{a},P,M,X))\biggl\}\\
\end{array}
$$
\end{lemma}

\noindent
(Proof). Let $\bar{a}$ be as stated. Suppose that 
$$
conf(\bar{a}_0)=C_{INIT}(M,X)\land
\all i<s\delta_M(P(i),conf(\bar{a}_i),conf\bar{a}_{i+1}))\land
Accept(conf(\bar{a}_s),M,X)).
$$
By induction on $k$ we show that 
$\all k<l_0sg((G_P)_{merge(\bar{a}^{<l_0-k},\bar{b}^{<l_0-k})})\neq 0$.
If $k=0$ then the claim follows from Corollary \ref{cor:accepting2-1} since 
$$
sg(((G_P)_{merge(\bar{a}^{<l_0},\bar{b}^{<l_0})})_{z^{l_0}})=sg((G_P)_{merge(\bar{a},\bar{b})})=0.
$$

For $k<l_0-1$, suppose by the inductive hypothesis that 
$sg((G_P)_{merge(\bar{a}^{<l_0-k},\bar{b}^{<l_0-k})})\neq 0$. 
Then 
$$
sg(((G_P)_{merge(\bar{a}^{<l_0-k},\bar{b}^{<l_0-k})})_{b_{l_0-k-1}})
=sg((G_P)_{merge(\bar{a}^{<l_0-k},\bar{b}^{<l_0-k})})\neq 0. 
$$
Thus by Corollary \ref{cor:legitimate2}, we have
$sg((G_P)_{merge(\bar{a}^{<l_0-k},\bar{b}^{<l_0-k-1})})\neq 0$. 
Since 
$$
sg(((G_P)_{merge(\bar{a}^{<l_0-k-1},\bar{b}^{<l_0-k-1})})_{a_{l_0-k-1}})
=sg((G_P)_{merge(\bar{a}^{<l_0-k},\bar{b}^{<l_0-k-1})}). 
$$
we have
$sg((G_P)_{merge(\bar{a}^{<l_0-(k+1)},\bar{b}^{<l_0-(k+1)})})\neq 0$ as desired. 

The converse direction is an immediate consequence of Corollary \ref{cor:accepting1}  and 
Corollary \ref{cor:accepting2-1}.\hfill$\Box$


\vspace{6pt}

Analogously, Alice always chooses options of $G_R(M,X)_P$ whose Sprague-Grundy number is 
equal to $0$ if and only if Bob's moves form a rejecting computation along $P$.  

\begin{lemma}\label{lemma:rejecting}
$\tpg$ proves that 
$$
\begin{array}{l}
\all M,X,P\all\bar{a}=\pairs{a_{0,0},\ldots a_{s,s+1}}
\biggl\{(|P|=s\land Leg(\bar{a}))\rightarrow\\
\all k<Len(\bar{a})(sg(G_R(M,X)_{P*merge(\bar{a}^{\leq k},\bar{b}^{<k})})=0)\leftrightarrow
RComp(\bar{a},P,M,X))\biggl\}\\
\end{array}
$$
\end{lemma}

\noindent
(Proof). Suppose that $RComp(\bar{b},P,M,X)$ holds. By induction on $k$, 
we show that $\all k<(s+1)(s+2)sg(G(M,X)_{P*merge(\bar{a}^{\leq (s+1)(s+2)-k},\bar{b}^{(s+1)(s+2)<k})})\neq 0$. 
If $k=0$ then the claim follows from Corollary \ref{cor:accepting2-2} since
$$
sg(G(M,X)_{P*merge(\bar{a}^{\leq (s+1)(s+2)-k},\bar{b}^{(s+1)(s+2)<k})*b_{s,i}^Q})=0
$$
for $i\neq 1$. The proof for $k>0$ is identical to the one for Lemma \ref{lemma:accepting}. 
\hfill$\Box$

Finally we show that applying the strategy function $\tau_A$ to either $G_A(M,X)_P$ or 
$G_R(M,X)_P$ yields either accepting or rejecting computation respectively. 

\begin{lemma}\label{lemma:accepting2}
$\tpg$ proves that if $sg(G_A(M,X)_P)\neq 0$ then the application of $\tau_A$ to 
$G_A(M,X)_P$ yields an accepting computation along $P$: 
$$
\begin{array}{l}
\all M,X,P\all\bar{a}=\pairs{a_{0,0},\ldots a_{s,s+1}}\\
\left\{(|P|=s\land sg(G_A(M,X)_P)\neq 0\land
\tau_A(\bar{b},G_A(M,X)_P)=merge(\bar{a},\bar{b}))
\rightarrow AComp(\bar{a},P,M,X)\right\}
\end{array}
$$
\end{lemma}

\noindent
(Proof). Suppose that $sg(G_A(M,X)_P)\neq 0$ and let 
$\tau_A(\bar{b},G_A(M,X)_P)=merge(\bar{a},\bar{b})$. By the definition of $\tau_A$, we have 
$$
\all k\leq Len(\bar{a})(sg((G_P)_{merge(\bar{a}_{\leq k},\bar{b}_{<k})})=0).
$$
So by Lemma \ref{lemma:accepting}, we have the claim. \hfill$\Box$


\begin{lemma}\label{lemma:rejecting2}
$\tpg$ proves that if $sg(G_R(M,X)_P)\neq 0$ then the application of 
$\tau_A$ to $G_R(M,X)_P$ yields a rejecting computation: 
$$
\begin{array}{l}
\all M,X,P\all\bar{a}=\pairs{a_{0,0},\ldots a_{s,s+1}}\\
\left\{(|P|=s\land sg(G_R(M,X)_P)\neq 0\land\tau_A(\bar{a},M,X)=merge(\bar{a},\bar{b}))
\rightarrow RComp(\bar{a},P,M,X)\right\}
\end{array}
$$
\end{lemma}

\noindent
(Proof). Suppose that $sg(G(M,X)_P\neq 0$. Then By Lemma \ref{lemma:rejecting} 
we have the claim by a similar argument as for Lemma \ref{lemma:accepting2}. \hfill$\Box$


\vspace{6pt}

Next lemma states that $G_A(M,X)_P$ and $G_R(M,X)_P$ play complementary roles to each other. 

\begin{lemma}\label{lemma:complement}
$\tpg$ proves that 
$$
\all M,X,P(|P|=s\rightarrow
(sg(G_A(M,X)_P)\neq 0\leftrightarrow sg(G_A(M,X)_P)=0)).
$$
\end{lemma}

\noindent
(Proof). We argue in $\tpg$. Suppose that $sg(G_A(M,X)_P)\neq 0$. We show that 
$$
\begin{array}{l}
\all\bar{a}(Len(\bar{a})=(s+1)(s+2)\rightarrow\\
\ext\bar{b}(Len(\bar{b})\leq Len(\bar{a})\land
WSeq(merge(\bar{a}^{\leq Len(\bar{b})},\bar{b}),G_R(M,X)_P))).
\end{array}(*)
$$
The proof is divided into cases. Let $\bar{a}$ be an arbitrary list of Alice's moves with 
$Len(\bar{a})=(s+1)(s+2)$ and $\bar{b'}=\pairs{b_{0,0},\ldots,b_{s,s}}$. 

If $A\mbox{-}Leg(\bar{a})\land Comp(\bar{a})$ then by Corollary \ref{cor:accepting2-1}, we have 
$$
G_A(M,X)_{P*merge(\bar{a},\bar{b'})}=\emptyset\leftrightarrow AComp(\bar{a},M,X,P). 
$$
On the other hand, by Corollary \ref{cor:accepting2-2}, we have 
$$
G_R(M,X)_{P*merge(\bar{a},\bar{b'})}=\emptyset\leftrightarrow RComp(\bar{a},M,X,P). 
$$
Thus we have $G_R(M,X)_{P*merge(\bar{a},\bar{b'})}\neq\emptyset$ and for any 
$c\in Node(G_R(M,X)_{P*merge(\bar{a},\bar{b'})}\subseteq C_R$, we have 
$G_R(M,X)_{P*merge(\bar{a},\bar{b'}*c)}=\emptyset$. 
Therefore we obtain $WSeq(merge(\bar{a},\bar{b}*c),G_R(M,X)_P)))$. 

If $A\mbox{-}Leg(\bar{a})\land\neg Comp(\bar{a})$ then by Corollary \ref{cor:accepting1} 
we have 
$$
G_R(M,X)_{P*merge(\bar{a},\bar{b'})}\neq\emptyset\land
G_R(M,X)_{P*merge(\bar{a},\bar{b'})*c}=\emptyset.
$$

Finally if $\neg A\mbox{-}Leg(a)$ then we can find the shortest initial part 
$\bar{a'}=\pairs{a_0,\ldots,a_k}$ of $\bar{a}$ such that 
$A\mbox{-}Leg(\bar{a'})\land a_{k+1}\not\in A_{q_{k+1},r_{k+1}}$. 
Then by Lemma \ref{lemma:legitimate1}, we have $x$ such that 
$$
WSeq(merge(\bar{a'},\bar{b}^{\leq k})*a_{k+1}*x,G_R(M,X)_P).
$$

Thus in any case we have $(*)$ and from this we readily have $sg(G_R(M,X))=0$. 
Conversely, if $sg(G_R(M,X))\neq 0$ the by a similar argument, we obtain $sg(G_A(M,X))=0$. 
\hfill$\Box$

\vspace{6pt}

\noindent
(Proof of Lemma \ref{lemma:comp2}). Suppose that $sg(G_A(M,X))\neq 0$. 
Then by Lemma \ref{lemma:accepting2}, we have the first part. 
If $sg(G_A(M,X))=0$ then by Lemma \ref{lemma:complement}, we have $sg(G_R(M,X))\neq 0$ and 
we can apply Lemma \ref{lemma:rejecting2} 

\vspace{6pt}
\noindent
(Proof of Theorem \ref{theorem:main}). We argue in $\tpg$. 
Let $M$ be an alternating Turing machine and $X$ be an input. We define $G(M,X)=G_A(M,X)$ 
For other functions, we set 
$$
\begin{array}{l}
Path_A(P,M,X)=\tau_A(P,G_A(M,X))\\
Path_R(P,M,X)=\tau_B(P,G_A(M,X))\\
Comp_A(M,X,P)=ASeq(\tau_A(\bar{b'},G_A(M,X)_P))
Comp_R(M,X,P)=ASeq(\tau_A(\bar{b'},G_A(M,X)_P))
\end{array}
$$
where the sequence $\bar{b}$ is defined by $\bar{b}=\pairs{b_{0,0},\ldots,b_{s,s}}$. 

The condition (1) is trivial from the definition. Conditions (2) and (3) follows from 
the definition of the strategy functions $\tau_A$ and $\tau_B$. 

Since we assume that $p(|X|)$ is even for all $X$, it follows that
$$
\all X,P(|P|=p(|X|)\rightarrow(sg(G(M,X)=0\leftrightarrow sg(G_A(M,X)_P)=0)).
$$
Thus Lemma \ref{lemma:comp2} implies 4. So the proof terminates.\hfill$\Box$

\begin{theorem}
$\tpg$ proves $\sbi{\infty}$-IND. 
\end{theorem}

\noindent
(Proof). 
For any $\varphi(X)\in\sbi{\infty}$ we can construct an alternating Turing machine 
which decides $\varphi$ in polynomial time. \hfill$\Box$

\end{document}